\input amstex
\documentstyle{amsppt}
\magnification=\magstep1 \NoRunningHeads

\topmatter

\title
 Rank-one flows of transformations with infinite ergodic index
 \endtitle
\author  Alexandre~I.~Danilenko and Kyewon K. Park
\endauthor
\abstract
A rank-one infinite measure preserving  flow $T=(T_t)_{t\in\Bbb R}$ is constructed  such  that for each $t\ne 0$, the Cartesian powers of the transformation $T_t$ are all ergodic.
\endabstract

\address
 Institute for Low Temperature Physics
\& Engineering of National Academy of Sciences of Ukraine, 47 Lenin Ave.,
 Kharkov, 61164, UKRAINE
\endaddress
\email alexandre.danilenko\@gmail.com
\endemail

\address Department of Mathematics, College of Natural Science,
Ajou University, Suwon 442-749, KOREA
\endaddress
\email kkpark\@{madang.ajou.ac.kr}
\endemail

\keywords  Infinite ergodic index, $(C,F)$-construction
 \endkeywords

\subjclass Primary 37A40; Secondary 37A10, 37A25
\endsubjclass
\thanks
The first author thanks Ajou University and KIAS for  supporting in part his visit to  South Korea.
The second author is supported in part by KRF 2007-313-C00044.
\endthanks

\endtopmatter
\document

\head 0. Introduction
\endhead

In  1963 Kakutani and Parry  discovered an interesting phenomenon in the theory of infinite measure preserving  maps.
They showed that for each $p>0$,  there exists a transformation whose $p$-th Cartesian power is ergodic but $(p+1)$-th one is not \cite{KP}.
Since then a number of other examples of transformations with exotic (from the point of view of the classical ``probability preserving'' ergodic theory)
weak mixing properties were constructed.
 See surveys \cite{Da2} and \cite{DaS3} for a detail discussion on that.
In \cite{Da1}, \cite{DaS1} these examples were extended to infinite measure preserving actions of discrete countable Abelian groups.
Weak mixing properties of infinite measure preserving actions of {\it continuous} Abelian groups such as $\Bbb R$ and $\Bbb R^d$ were under consideration in \cite{I--W}.
In particular, a rank-one flow (i.e. $\Bbb R$-action) whose Cartesian square is ergodic was constructed there. A rank-one infinite measure preserving flow $T=(T_t)_{t\in\Bbb R}$ with infinite ergodic index (i.e. the Cartesian powers of $T$ are  all ergodic) appeared in a recent paper \cite{DaSo}.
It can be deduced easily  from \cite{DaSo} that there is a residual subset $D_T$ of $\Bbb R$ such that  for each $t\in D_T$, the transformation $T_t$ has infinite ergodic index.  However the following more subtle question by C.~Silva remained open so far:
\roster
{\it Is there a rank-one infinite measure preserving flow $T$ with $D_T=\Bbb R\setminus\{0\}$?}
\endroster
Our purpose in this paper is to answer his question in the affirmative.

\proclaim{Theorem 0.1}
There is a rank-one infinite $\sigma$-finite measure preserving flow $T=(T_t)_{t\in\Bbb R}$ such that for each $t\ne 0$, the transformation $T_t$ has infinite ergodic index.
\endproclaim

The main idea of the proof is different from those that were used in \cite{I--W} and \cite{DaSo}. It is based on a technique to {\it force} a dynamical property.
Originating from \cite{Ry1}, such  techniques were  utilized in \cite{Ry2}, \cite{DaR}, etc. to obtain mixing, power weak mixing, etc. of some systems.
In this paper the desired flow appears as a certain limit of a sequence of weakly mixing finite measure preserving flows.
We construct this sequence in such a way to retain the property of infinite ergodic index in the limit.
The construction is implemented in the language of $(C,F)$-actions (see \cite{Da2}).

\head 1. Preliminaries: rank-one actions and $(C,F)$-actions of $\Bbb R^d$
\endhead
We first recall the definition of rank one.
Let $S=(S_g)_{g\in \Bbb R^d}$ be a measure preserving action of $\Bbb R^d$ on a standard
$\sigma$-finite measure space $(Y,\goth C,\nu)$.

\definition{Definition 1.1}
\roster
\item"(i)"
A {\it Rokhlin tower or column} for $S$ is a triple $(A,f, F)$,
 where $A\in\goth C$, $F$ is a cube in $\Bbb R^d$
and $f:A\to F$ is a measurable mapping
 such that for any Borel subset $H\subset F$ and an element $g\in \Bbb R^d$ with
 $g+H\subset F$, one has
 $f^{-1}(g+H)=S_gf^{-1}(H)$.
\item"(ii)"
 $S$ is said to be of {\it rank-one (by cubes)} if
there exists a sequence of Rokhlin towers $(A_n,f_n,F_n)$ such that
the volume of $F_n$ goes to infinity and
for any subset $C\in\goth C$ of finite measure, there is a sequence
of Borel subsets $H_n\subset F_n$ such that
$$
\lim_{n\to\infty}\nu(C\triangle f_n^{-1}(H_n))=0.
$$
\endroster
\enddefinition

The $(C,F)$-construction of measure preserving actions for discrete countable groups was introduced in \cite{dJ} and \cite{Da1}. It was extended to the case of locally compact second countable Abelian groups in \cite{DaS2}. (See also  \cite{Da2}.) Here we outline it briefly for $\Bbb R^d$, $d\in\Bbb N$.

Given two subsets $E,F\subset \Bbb R^d$,  by $E+F$ we mean their algebraic sum, i.e. $E+F=\{\,e+f\mid e\in E,f\in F\,\}$.
The algebraic difference $E-F$ is defined in a similar way.
If $F$ is a singleton, say $F=\{f\}$, then we will write $E+f$ for $E+F$.
Two subsets $E$ and $F$ of $\Bbb R^d$ are called {\it independent} if $(E-E)\cap(F-F)=\{0\}$, i.e. if $e+f=e'+f'$ for some $e,e'\in E$, $f,f'\in F$ then $e=e'$ and $f=f'$.

Fix $p\in\Bbb N$ and consider two sequences $(F_n)_{n=0}^\infty$ and $(C_n)_{n=1}^\infty$ of subsets in $\Bbb R^d$ such that $F_n$ is a  cube  $[0,h_n)\times\cdots\times[0,h_n)$ ($d$ times) for an $h_n\in\Bbb R$, $C_n\subset\Bbb R^d$ is a finite set, $\#C_n > 1$,
$$
\align
&F_n \text{ and } C_{n+1} \text{ are independent};\tag1-1\\
&F_n + C_{n+1}\subset F_{n+1}.\tag1-2
\endalign
$$
This means that $F_n + C_{n+1}$ consists of $\#C_{n+1}$ mutually disjoint `copies' $F_n+c$ of $F_n$, $c\in C_{n+1}$,  and all these copies are contained in $F_{n+1}$.
We equip $F_n$ with the measure $(\#C_1\cdots\# C_n)^{-1}(\lambda_{\Bbb  R^d}\restriction F_n)$, where $\lambda_{\Bbb R^d}$ denotes Lebesgue measure on $\Bbb R^d$. Endow $C_n$ with the equidistributed probability measure.
Let $X_n := F_n\times\prod_{k>n}C_k$ stand for the product of measure spaces.
Define an embedding $X_n\to X_{n+1}$ by setting
$$
(f_n, c_{n+1}, c_{n+2}, \ldots )\mapsto(f_n+c_{n+1}, c_{n+2}, \ldots ).
$$
 It is easy to see that this embedding is measure preserving.
Then $X_0 \subset X_1\subset\cdots$. Let $X:=\bigcup_{n=0}^\infty X_n$
denote the inductive limit of the sequence of measure spaces $X_n$ and let
$\mu$ denote the corresponding measure on $X$. Then $\mu$ is
$\sigma$-finite. It is
 infinite if and only if
$$
\lim_{n\rightarrow\infty}\frac{h_n^d}{\#C_1\cdots\#C_n}=\infty.
\tag1-3
$$
Given $g\in\Bbb R^d$ and $n\in\Bbb N$, we set
$$
L_g^{(n)}:=(F_n\cap(F_n-g))\times\prod_{k>n}C_k\quad\text{and}\quad
R_g^{(n)}:=(F_n\cap(F_n+g))\times\prod_{k>n}C_k.
$$
Clearly, $L_g^{(n)}\subset L_g^{(n+1)}$ and $R_g^{(n)}\subset R_g^{(n+1)}$.
Define a map $T_g^{(n)}\colon L_g^{(n)}\to R_g^{(n)}$ by setting
$$
T_g^{(n)}(f_n,c_{n+1},\ldots):=(f_n+g,c_{n+1},\ldots).
$$
Put
$$
L_g:=\bigcup_{n=1}^\infty L_g^{(n)}\subset X\quad\text{and}\quad
R_g:=\bigcup_{n=1}^\infty R_g^{(n)}\subset X.
$$
Then a Borel one-to-one map $T_g\colon L_g\to R_g$ is well defined by $T_g\upharpoonright L_g^{(n)}=T_g^{(n)}$.
Since $h_n\to\infty$, it follows that $\mu(X\setminus L_g)=\mu(X\setminus R_g)=X$ for each $g\in\Bbb R^d$.
It is easy to verify that  $T:=(T_g)_{g\in\Bbb R^d}$ is a Borel $\mu$-preserving action of $\Bbb R^d$.

\definition{Definition 1.2}
$T$ is called the $(C,F)$-{\it action of $\Bbb R^d$ associated with} $(C_{n+1}, F_{n})_{n\ge 0}$.
\enddefinition

Each $(C,F)$-action is of  rank one.

Given a Borel subset $A\subset F_n$, we set $[A]_n:=\{\,x=(x_i)_{i=n}^\infty\in X_n\mid x_n\in A\,\}$ and call it an {\it $n$-cylinder} in $X$.
Clearly,
$$
[A]_n=\bigsqcup_{c\in C_{n+1}}[A+c]_{n+1}.
$$
Notice also that
$$
T_g[A]_n=[A+g]_n \text{ for all }g\in\Bbb R^d\text{ and }A\subset F_n\cap(F_n-g), \ n\in\Bbb N.
\tag1-5
$$
The  sequence  of all $n$-cylinders approximates the entire Borel $\sigma$-algebra on $X$ when $n\to\infty$.

We state without proof the following standard lemma (see, e.g., Lemma~2.4 from \cite{Da1}).

\proclaim{Lemma 1.3} Let $\Cal P_n$ be a finite partition of $F_n$ into parallelepipeds   such that
for each atom $\Delta$ of $\Cal P_n$ and an element $c\in C_{n+1}$,
the parallelepiped $\Delta+c$ is $\Cal P_{n+1}$-measurable and the maximal diameter of the atoms in $\Cal P_n$ goes to zero.
Let $S$ be a measure preserving transformation of $X$.
Then the following holds.
\roster
\item"\rom{(i)}"
The sequence of collections of $n$-cylinders $\{[A]_n\mid A\subset F_n\text{ is $\Cal P_n$-measurable}\}$
approximates  the entire $\sigma$-algebra $\goth B$ as $n\to\infty$.
\item"\rom{(ii)}"
If for each pair of atoms $\Delta_1,\Delta_2\in\Cal P_n$, there are a subset $A\subset[\Delta_1]_n$  and a $\mu$-preserving one-to-one map $\gamma:A\to [\Delta_2]_n$ such that
$\mu(A)>0.5\mu([\Delta_1]_n)$,
and $\gamma x\in\{S^ix\mid i\in\Bbb Z\}$ for all $x\in A$
then $S$ is ergodic.
\endroster
\endproclaim

We will also use the following  property of the $(C,F)$-actions.
If $T$ is associated with $(C_{n+1}, F_n)_{n\ge 0}$ then for each $p>1$, the product action
$$
(T_{t_1}\times\cdots\times T_{t_p})_{(t_1,\dots,t_p)\in(\Bbb R^d)^p}
$$
is the $(C,F)$-action of $(\Bbb R^d)^p$ associated with  $(C_{n+1}^p,F_n^p)_{n\ge 0}$. The upper index $p$ means the $p$-th Cartesian power.

\head
2. Two auxiliary facts
\endhead

Given a $\sigma$-finite measure space $(X,\mu)$, we denote by Aut$(X,\mu)$ the group of all $\mu$-preserving (invertible)
transformations of $X$.
It is a Polish group when endowed with the {\it weak topology} \cite{Aa}.
Recall that the weak topology is the weakest topology in which the maps
$$
\text{Aut}(X,\mu)\ni T\mapsto\mu(TA\cap B)\in\Bbb R
  $$
are continuous for all subsets $A,B\subset X$ of finite measure.

 Given   $S\in\text{Aut}(X,\mu)$  and two subsets $A,B\subset X$ with $\mu(A)=\mu(B)<\infty$, we define  subsets $A_0,A_1,\dots$  of $A$ as follows:
$$
\aligned
A_0&:=A\cap B,\\
A_i&:=\bigg(A\setminus\bigsqcup_{j=0}^{i-1}A_j\bigg)\cap S^{-i}\bigg(B\setminus \bigsqcup_{j=0}^{i-1}S^jA_j\bigg), \quad i>0.
\endaligned
\tag2-1
$$
We now let $\Cal N_{S,A,B}:=\min\{i\ge 0\mid \mu(A_0\sqcup\dots\sqcup A_i)>0.5\mu(A)\}$.
If $S$ is ergodic then $A=\bigsqcup_{i\ge 0}A_i$ and hence $\Cal N_{S,A,B}$ is well defined.
Denote by $\Cal E$ the subset of all ergodic transformations  in
$
\text{Aut}(X,\mu)$.
It is well known that $\Cal E$ is a dense $G_\delta$ in $\text{Aut}(X,\mu)$.
Since for each $i\ge 0$, the map
$$
\text{Aut}(X,\mu)\ni S\mapsto\mu(A_0\sqcup\dots \sqcup A_i)\in\Bbb R
$$
is continuous, we obtain the following lemma.

\proclaim{Lemma 2.1}
The map
$
\Cal E\ni S\mapsto \Cal N_{S,A,B}\in\Bbb R
$
is upper semicontinuous for all subsets $A,B\subset X$ with $\mu(A)=\mu(B)<\infty$.
\endproclaim

In the case of $(C,F)$-actions we can say more about the ``structure''
of the sets $A_i$, $i=0,\dots,\Cal N_{S,A,B}$. For $q=(q_1,\dots,q_d)\in\Bbb R^d$, we let $\|q\|:=\max_{1\le i\le d}|q_i|$.

\proclaim{Lemma 2.2} Let \, $(X,\mu,(T_t)_{t\in\Bbb R^d})$ be a $(C,F)$-action of $\Bbb R^d$ associated with a sequence $(C_{n+1},F_n)_{n\ge 0}$ such that
$$
a+F_n+C_{n+1}\subset F_{n+1}\tag2-2
$$
for each $a=(a_1,\dots,a_p)$ with $a_i\ge 0$, $i=1,\dots,p$, and $\|a\|\le 1$.
Fix two $n$-cylinders
 $A$ and $B$ of equal measure  and a transformation $S=T_q$ for some $q\in\Bbb R^d_+$.
Then the subsets $A_0,A_1,\dots, A_{\Cal N(S,A,B)}$ defined by \thetag{2-1} are  $(n+Q\cdot\Cal N(S,A,B)
)$-cylinders, where $Q$ is any integer greater than $\|q\|$.
\endproclaim

\demo{Proof}
We let $N:=Q\cdot\Cal N(S,A,B)$.
If $A=[\widetilde A]_n$ for some $\widetilde A\subset F_n$ then
$A=[\widehat A]_{n+N}$, where $\widehat A:=\widetilde A+C_{n+1}+\cdots+C_{n+N}\subset F_{n+N}$.
From \thetag{1-2} and \thetag{2-2} we deduce that the sets $\widehat A+q,\dots,\widehat A+\Cal N(S,A,B)q$ are all contained in $F_{n+N}$.
It remains to use~\thetag{2-1} and \thetag{1-5}.  \qed
\enddemo

We also note that $S^iA_i\subset B$ and $S^iA_i\cap S^jA_j=\emptyset$ for all $i,j=0,\dots,\Cal N(S,A,B)$.

\head 3. Proof of the main result
\endhead

\proclaim{Theorem 3.1} There exists a $(C,F)$-flow $T=(T_t)_{t\in\Bbb R}$ such that each transformation $T_t$, $t\ne 0$, has infinite ergodic index.
\endproclaim
\demo{Proof}
We will construct this flow via an inductive procedure.
Fix a sequence of integers $(p_n)_{n\ge 1}$ in which every integer greater then 1
occurs infinitely many times.
Suppose that after   $n-1$  steps  of the construction we have already defined
$$
F_0,C_1,F_1,\dots, C_{m_{n-1}}, F_{m_{n-1}}.\tag3-1
$$
Suppose also that for each $0\le i\le m_{n-1}$, a finite partition $\Cal P_i$ of
$F_i$ into intervals  is chosen in such a way that
\roster
\item"---" the interval $\Delta+c$ is $\Cal P_{i+1}$-measurable for each atom $\Delta$ of $\Cal P_i$, $0\le i<m_{n-1}$, and each $c\in C_{i+1}$ and
\item"---" the length of any atom of $\Cal P_i$ is no more than $i^{-1}$, $1\le i\le m_{n-1}$.
\endroster
{\it Step $n$}. Consider a  rank-one weakly mixing finite measure preserving $(C,F)$-flow $T^{(n)}=(T_{t}^{(n)})_{t\in\Bbb R}$ associated with a sequence $(C_{k+1,n},F_{k,n})_{k\ge 0}$
such that $F_{0,n}:=F_{m_{n-1}}$.
 Examples of weakly mixing rank-one finite measure preserving flows are well known---see, e.g., \cite{dJP}.
 In  \cite{DaS2} one can find explicit  $(C,F)$-construction of mixing
finite measure preserving flows.
Let $(X^{(n)},\mu_{n})$ be the space of this action.
Since  $T^{(n)}$ is  weakly mixing, it follows that for each $t>0$, the  transformation
$$
S_t:=T_{t}^{(n)}\times\cdots\times T_{t}^{(n)} \text{ ($p_{n}$ times)}
$$
of the product space $(X^{(n)},\mu_{n})^{p_{n}}$ is ergodic.
We note that this space is the space of the $(C,F)$-action of $\Bbb R^{p_{n}}$ associated with the sequence $(C_{k+1,n}^{p_{n}},F_{k,n}^{p_{n}})_{k\ge 0}$
(see our remark at the end of \S1).
Given $k\ge 0$, let $\Cal P_{k,n}$ be a finite partition of $F_{k,n}$
into intervals  such that
\roster
\item"---" $\Cal P_{0,n}=\Cal P_{m_{n-1}}$,
\item"---" the interval $\Delta+c$ is $\Cal P_{k+1,n}$-measurable for each atom $\Delta$ of $\Cal P_{k,n}$ and  $c\in C_{k+1,n}$ and
\item"---" the length of any atom of $\Cal P_{k,n}$ is no more than $(m_{n-1}+k)^{-1}$.
\endroster
We now let $\Cal P_{k,n}^{p_{n}}:=\Cal P_{k,n}\times\cdots\times \Cal P_{k,n}$ ($p_{n}$ times).
Then $\Cal P_{k,n}^{p_{n}}$ is a finite partition of $F_{k,n}^{p_{n}}$ into parallelepipeds and
\roster
\item"---" the parallelepiped $\Delta+c$ is $\Cal P_{k+1,n}^{p_{n}}$-measurable for each atom $\Delta$ of $\Cal P_{k,n}^{p_{n}}$ and  $c\in C_{k+1,n}^{p_{n}}$ and
\item"---" the diameter of an atom of $\Cal P_{k,n}^{p_{n}}$ is no more than $(m_{n-1}+k)^{-p_{n}}$.
\endroster
Denote by $D_n$ the maximum of $\Cal N(S_t,[\Delta]_0,[\Delta']_0)$ when $\Delta$ and $\Delta'$ run independently the atoms of
$\Cal P_{0,n}^{p_{n}}$
 and $t$ runs the segment $[n^{-1},n]\subset\Bbb R$.
It exists by Lemma~2.1.
It now follows from Lemma~2.2 that  for any
pair of parallelepipeds  $\Delta,\Delta'\in\Cal P_{0,n}^{p_{n}}$ and a real $t\in [n^{-1},n]$, there exist $nD_n$-cylinders
 $A_1,\dots, A_{D_n}\subset [\Delta]_0$
such that
$$
\gathered
\mu_{n}^{p_n}\bigg(\bigsqcup_{i=1}^{D_n}A_i\bigg)>\frac 12\mu_{n}^{p_n}([\Delta]_0),\\
S_t^iA_i\subset [\Delta']_0\text{ for each $1\le i\le D_n$ and} \\
S_t^iA_i\cap S_t^jA_j=\emptyset\text{ if }1\le i\ne j\le D_n.
\endgathered
\tag3-2
$$
We now ``continue'' the sequence \thetag{3-1} by setting
$$
C_{m_{n-1}+1}:=C_{1,n}, F_{m_{n-1}+1}:=F_{1,n},\dots, C_{m_{n-1}+nD_n}:=C_{nD_n,n}.
$$
Next, to define $F_{m_{n-1}+nD_n}$ we ``double'' the set  $F_{nD_n,n}$, i.e.
$$
F_{m_{n-1}+nD_n}:=[0,2a)\text{ if }F_{nD_n,n}=[0,a)\text{ for some }a>0.\tag3-3
$$
It remains to put $m_{n}:=m_{n-1}+nD_n$. The $n$-th step is now completed.

Continuing this procedure  infinitely many times, we obtain the entire sequence $(C_{i+1},F_{i})_{i=0}^\infty$.
Denote by $T=(T_t)_{t\in\Bbb R}$ the associated $(C,F)$-flow.
Let $(X,\mu)$ be the space of this flow.
It follows from \thetag{3-3} that
$\lambda_\Bbb R(F_i)>2\lambda_\Bbb R(F_{i-1})\# C_i$ for infinitely many $i$. Hence
$\mu(X)=\infty$.
Moreover, a finite partition $\Cal P_i$ of $F_i$ into intervals  is fixed such that the conditions of Lemma~1.3 are satisfied.
Next, there are  one-to-one  correspondences (natural identifications) between
\roster
\item"---"  the collection of 0-cylinders in $X^{(n)}$ and
the collection of $m_{n-1}$-cylinders in $X$ and
\item"---" the collection of $nD_n$-cylinders in $X^{(n)}$ and
the collection of $m_{n}$-cylinders in $X$.
\endroster
Moreover, the ``dynamics'' of $T^{(n)}$ on the $nD_n$-cylinders is the same as the dynamics of $T$ on the $m_{n}$-cylinders.
This means the following: if $A,B\subset F_{nD_n}^{(n)}$ and $[B]_{nD_n}=T_{w}^{(n)}[A]_{nD_n}$ for some $w\in\Bbb R$ then $[B]_{m_{n}}=T_{w}[A]_{m_{n}}$.
Therefore we deduce from \thetag{3-2} that
for any
pair of parallelepipeds $\Delta,\Delta'\in\Cal P_{m_n-1}^{p_{n}}$ and a real $t\in [n^{-1},n]$, there exist $m_{n}$-cylinders
 $A_1,\dots, A_{D_n}\subset [\Delta]_{m_{n-1}}$
such that
$$
\gathered
\mu^{p_n}\bigg(\bigsqcup_{i=1}^{D_n}A_i\bigg)>\frac 12\mu^{p_n}([\Delta]_{m_{n-1}}),\\
V_t^iA_i\subset [\Delta']_{m_{n-1}}\text{ for each $1\le i\le D_n$ and} \\
V_t^iA_i\cap V_t^jA_j=\emptyset\text{ if }1\le i\ne j\le D_n,
\endgathered
$$
where $V_t:=T_t\times\cdots\times T_t$ ($p_{n}$ times).
Fix $p>0$.
Passing to a subsequence where $p_{n}=p$ we now deduce
 from Lemma~1.3(ii) that $T_t\times\cdots\times T_t$ ($p$ times) is ergodic for each $t>0$.
Hence $T_t$ has infinite ergodic index.
\qed
\enddemo

\head 4. Concluding remarks
\endhead

\subhead 4.1\endsubhead When constructing $T$, we use only finitely many initial terms of the sequence $(C_{k+1}^{(n)},F_k^{(n)})_{k\ge 0}$ for each $n>0$.
However to determine ``where to stop'' (i.e. to determine $D_n$) we use the weak mixing properties of the auxiliary flow $T^{(n)}$ which depends on the entire {\it infinite} sequence $(C_{k+1}^{(n)}, F_k^{(n)})_{k\ge 0}$. No upper bound on $D_n$ is found.
This means that the construction of $T$ is not {\it effective}.
In this connection  we  rise a  question:
\roster
{\it is it possible to find an effective construction for the flow from Theorem~0.1?}
\endroster
We note that the construction in \cite{DaSo} is effective.

\subhead 4.2\endsubhead
It is possible to strengthen Theorem 0.1 by replacing  the infinite ergodic index with a stronger property of power ergodicity. Recall that a measure preserving transformation $S$ is called power ergodic if for each finite sequence $n_1,\dots,n_k$ of nonzero integers the transformation
$S^{n_1}\times\cdots\times S^{n_k}$ is ergodic.
Only  a slight modification of our argument is needed to   show the following theorem.

\proclaim{Theorem 4.1}
{\it There exists a rank-one infinite $\sigma$-finite measure preserving flow $T=(T_t)_{t\in\Bbb R}$  such that
the transformation $T_t$ is power weakly mixing for each $t\ne 0$.}
\endproclaim

Also, it is easy to extend Theorem 0.1 to actions of  $\Bbb R^d$.

\proclaim{Theorem 4.2}
For each $d>1$, there exists a rank-one infinite $\sigma$-finite measure preserving action $T=(T_g)_{g\in\Bbb R^d}$ of $\Bbb R^d$ such that
the transformation $T_g$ has infinite ergodic index for each $g\ne 0$.
\endproclaim

We leave the proofs of Theorems 4.1 and 4.2 to the reader.

\enddocument